\documentclass{article}
\usepackage{diagrams}

\usepackage{mathptmx}

\usepackage[textwidth=12.8cm, textheight=18.5cm]{geometry}

\usepackage{amssymb}
\usepackage{amsmath}

\newtheorem{prop}{Proposition}
\newtheorem{thm}{Theorem}
\newtheorem{corol}{Corollary}

\title{Calculus of extensive quantities}
\author{Anders Kock\\
University of Aarhus}
\date{}

\DeclareMathOperator{\tot}{tot}

\newcommand{\E}{\mathcal{E}}
\newcommand{\p}{\pitchfork}

\newcommand {\imes}{\times}
\begin{document}
\maketitle
\small \noindent {\bf Abstract.} We show how a commutative monad gives rise to a theory of 
extensive quantities, including (under suitable further conditions) a differential calculus of such.  The 
relationship to Schwartz distributions is dicussed. The paper is a 
companion to the author's ``Monads and extensive quantities'', but is 
phrased in more elementary terms.
\normalsize

\section*{Introduction}  Quantities of a given type distributed over a given space 
(say  distributions of smoke in a given room) 
may often be {\em added}, and multiplied by real scalars -- ideally, 
they form a real vector space. Lawvere 
stressed that the dependence of such vector spaces on the space over 
which the quantities in question are distributed, should be taken into account; in fact, 
the dependence 
is {\em functorial}. The viewpoint leads  to a distinction 
between two kinds of quantities: the functorality may be covariant, 
or it may be contravariant: In this context, the covariant quantity 
types are called {\em extensive quantities}, and the contravariant 
ones {\em intensive quantities}. This usage is an attempt to put 
mathematical precision into the  use of these terms in classical 
philosophy of physics. Mass distribution is an extensive quantity; 
mass density is an intensive one.
Lawvere observed that extensive and intensive quantites often come in 
pairs, with a definite pattern of mutual relationship, like the 
homology and cohomology functors on the category of topological 
spaces.

In \cite{MEQ}, we showed how  such a pattern essentially 
comes about, 
whenever one has a commutative monad $T$ on a Cartesian Closed Category 
$\E$ (where $\E$ is meant to model some category of {\em spaces}, not 
specified further). 

Such a monad is in particular a covariant endo-functor on $\E$, so 
the emphasis in our theory is the covariant aspect: the extensive 
quantities. We attempt to push  a theory of these as far as 
possible, with the intensive quantities in a secondary role.

This in particular applies to the {\em differential} calculus of extensive 
quantities on the line $R$, which will be discussed in the last 
Sections \ref{DCEQX} and \ref{EQSDX}; here, we also discuss the 
relationship to the  theory of Schwartz distributions of compact support 
(these have the covariant functorality requested for extensive 
quantities, and is a basis for  a classical version of a theory of extensive 
quantities). 

The article \cite{KRWE} by Reyes and the author develop some 
further differential calculus of extensive quantities, not only in 
dimension 1, as here; but it is couched in terms of the Schwartz 
(double dualization) paradigm, 
which we presently want to push in the background.

\section{Monads and their algebras}\label{MATAX}
The relationship between universal algebra, on the one hand, and 
monads\footnote{For the notion of monad, and algebras for a monad, 
the reader may consult \cite{Borceux}.} on 
the category of sets on the other, became apparent in the mid 60s, through 
the work of Linton, Manes, Kleisli, and many others:

If $T=(T, \eta , \mu )$ is such a monad, and $X$ is a set, an element 
$P\in T(X)$ may 
be interpreted as an $X$-ary operation on arbitrary $T$-algebras 
$B=(B, \beta )$: if $\phi :X\to B$ is an $X$-tuple of 
elements in $B$, we can construct a single element $\langle P, \phi 
\rangle \in B$, namely 
the value of 
$$\begin{diagram}T(X)&\rTo ^{T(\phi)}&T(B)&\rTo ^{\beta}& B
\end{diagram}$$
on the element $P\in T(X)$.
Then every morphism $f:B\to C$ of $T$-algebras is a homomorphism with 
respect to the operation defined by $P$. There is also a converse 
statement.

\begin{sloppypar}The monad-theoretic formulation of universal algebra can be lifted to 
 symmetric monoidal closed categories $\E$ other than sets, 
provided one considers the monad $T$ to be $\E$-enriched\footnote{for 
these notions, the reader is again referred to \cite{Borceux}.}, in 
particular, it applies to $\E$-enriched monads on any cartesian closed 
category $\E$.\end{sloppypar}

 Recall that for any functor $T:\E \to \E$, one has  maps $\hom 
(X,Y)\to \hom (T(X),T(Y))$, sending $f\in \hom (X,Y)$ to $T(f)\in \hom 
(T(X),T(Y))$; the $\E$-enrichment means that these maps not only can be 
defined for the $\hom${\em -sets}, but for the $\hom$-{\em objects}, 
so that we get  maps (the ``strength'' of $T$)
\begin{equation}\label{strx}st_{X,Y}: Y^{X}\to 
T(Y)^{T(X)};\end{equation}
 $Y^{X}$, as an object of $\E$ carries more structure than 
the mere {\em set} of maps from $X$ to $Y$, e.g.\ it may carry some topology, if 
$\E$ happens to be of topological nature.

This widening of the scope of monad theoretic universal algebra was documented 
in a series of articles by the author in the early 
1970s, cf.\ \cite{MSMCC}, \cite{DD}, \cite{CCGBCM}, \cite{BCCM}, 
\cite{SFMM}. In particular, the formulation makes sense for cartesian 
closed categories, 
which is 
the context of the present note. Via this formulation, it makes 
contact with functional analysis, because the basic logic behind 
functional analysis is the ability to form {\em function spaces}, 
even non-linear ones. For instance, the category $\E$ of convenient vector 
spaces, and the smooth (not necessarily non-linear maps) is a 
cartesian closed category, 
where several theories of functional analysis have  natural 
formulation, e.g.\ the theory of distributions (in the sense of 
Schwartz and others).

 We 
 expound here some aspects of the relationship between the theory of 
strong monads and functional analysis, by talking about $\E$ as if it 
were just 
the category of sets, and where $\E$-enrichment therefore is automatic. 
(A rigourous account of these aspects for general cartesian closed 
categories $\E$ is 
given in \cite{MEQ}.) So our exposition technique here is in the spirit of 
the ``naive'' exposition of synthetic differential geometry, as given in 
\cite{SGM}, say.

One aim of the theory developed in \cite{MEQ}, and expounded 
``synthetically'' in the present note, is to document that the 
space $T(X)$ may be seen as a space of {\em extensive 
quantities} (of some type) on $X$, in the sense of Lawvere. So we prefer to talk 
about a
$P \in T(X)$ as an {\em extensive quantity}  on 
$X$, rather than as ``an $X$-ary operation, operating on all 
$T$-algebras $B$''. To make this work well, one must assume that the 
monad $T$ is {\em commutative}\footnote{It should be stressed that to be commutative
is a property of {\em enriched} (=strong) monads, and enrichment is a {\em 
structure}, not a property. However, for the case where $\E$ is the 
category of sets, enrichment is automatic.}, see  Section \ref{TPEQX} and \ref{CMX}.

A main emphasis in Lawvere's ideas is that the space of extensive 
quantities on $X$ depend (covariant) {\em  functorially} on $X$; 
in the present context, functorality is encoded by 
the fact that $T$ is a functor. So for $f:X\to Y$, we have 
$T(f):T(X)\to T(Y)$; when $T$ is well understood from the context, we 
may write $f_{*}$ for $T(f)$; this is a type of notation that is as 
old as the very notion of functor (recall homology!). 

Thus, the semantics of $P\in T(X)$, as an $X$-ary operation an 
$T$-algebras $(B,\beta)$, may be rendered in terms of a pairing
\begin{equation}\label{semantx}\langle P, \phi \rangle := \beta 
(\phi_{*}(P)),\end{equation}
where  $\phi :X\to B$ is a map (``an $X$-tuple of elements of $B$''). This semantic aspect 
of extensive quantities is  essential in 
Section \ref{EQSDX}, where it is seen as the basis of a synthetic theory of Schwartz 
distributions.

For completeness, let us indicate how the pairing is defined without 
using individual {\em elements} $P\in T(X)$ and $\phi \in B^{X}$, as a 
map
$$T(X)\imes B^{X}\to B,$$
namely utilizing the assumed enrichment (\ref{strx}) of $T$  over $\E$:
$$\begin{diagram}T(X)\imes B^{X}&\rTo^{T(X)\times st}&T(X)\times 
T(B)^{T(X)}&\rTo ^{ev}&T(B)&\rTo^{\beta} &B;
\end{diagram}$$
here,  $ev$ denotes the ``evaluation'' map (part of the cartesian 
closed structure of $\E$). Henceforth, we shall be content with using 
``synthetic'' descriptions, utilizing elements.

We note the following naturality property of the pairing: for $f:X\to 
Y$, $P\in T(X)$ and $\psi :Y\to B$, we have
\begin{equation}\label{natux}
\langle f_{*}(P), \psi \rangle = \langle P, f^{*}(\psi )\rangle ,
\end{equation}
where $f^{*}(\psi ):= \psi \circ f$. For, the left hand side is 
$\beta (\psi_{*}f_{*}P)$, and the right hand side is 
$$\beta ((f^{*}(\psi ))_{*}P)= \beta ((\psi \circ f)_{*}(P)),$$
but $(\psi \circ f)_{*} = \psi _{*}\circ f_{*}$ since $T$ is a 
functor.

\begin{sloppypar}Similarly, if $B=(B,\beta )$ and $(C=(C,\gamma )$ are $T$-algebras, 
and $F:B\to C$ is a $T$-homomorphism\footnote{later in this article, 
$T$-homomorphisms will be called  ``{\em $T$-linear}'' maps.}, we have, for $P\in 
T(X)$ and $\phi : X\to B$ that
\begin{equation}\label{linx}
F(\langle P, \phi \rangle ) = \langle P , F\circ \phi \rangle ).
\end{equation}
This is an immediate consequence of $F\circ \beta =\gamma \circ T(F) 
$, the equation expressing that $F$ is a 
$T$-homomorphism.\end{sloppypar}

\medskip

The terminal object of $\E$ is denoted ${\bf 1}$. The object $T({\bf 1})$  
plays a special role as the algebra of ``scalars'', and we denote it 
also $R$; with suitable 
{\em properties} of $T$, it will in fact carry a canonical commutative ring 
structure, see Section \ref{ASX}. Its multiplicative unit $1$ is the 
element picked out by $\eta_{{\bf 1}}$. 
For any $X\in \E$, we 
have a unique map $X\to {\bf 1}$, denoted $!$ (when $X$ is understood from 
the context). For $P\in T(X)$, we have canonically associated a 
scalar $\tot (P)\in T({\bf 1})$, the {\em total} of $P$, namely
$$\tot (P):= !_{*}(P).$$ From uniqueness of maps to ${\bf 1}$ follows 
immediately that $P$ and $f_{*}(P)$ have the same total, for any 
$f:X\to Y$.

One example of a monad  $T$ (on a suitable cartesian closed category 
of smooth spaces) is where $T(X)$ is the space of Schwartz 
distributions {\em of compact support}; we return to Schwartz 
distributions in \ref{EQSDX}, and they are only mentioned here as a warning, 
namely that  functorality is a strong requirement; thus for 
instance, a uniform  distribution on the line can never have a total.
 (The functorial properties of non-compact distributions are not  
understood well enough presently.) In \cite{BET}, the authors 
construct the free real vector space monad $T$, in a category of suitable 
``convenient'' spaces (\cite{KM}), by carving it, out by topological 
means, from the monad of  
(compactly supported) Schwartz distributions. 

The units $\eta _{X} 
:X\to T(X)$ will, in terms of Schwartz distribution theory, pick out 
the Dirac distributions $\delta_{x}$; therefore, we shall allows 
ourselves the following doubling of notation: for $x\in X$, we write
$$ \eta _{X}(x) = \delta_{x}$$
(with $X$ understood from the context, on the right hand side).

\begin{prop}\label{pdelx}Let $B=(B,\beta )$ be a $T$-algebra, and 
$\phi:X\to B$ a map. Then
$$\langle \delta_{x},\phi \rangle = \phi (x).$$
In particular, 
$$\langle \delta_{x},\eta_{X} \rangle = \delta_{x}.$$  
\end{prop}
{\bf Proof.} In elementfree terms, the first equation  says that the composite
$$\begin{diagram}X&\rTo^{\eta_{X}}&T(X)&\rTo^{T(\phi)}&T(B)&\rTo 
^{\beta}&B
\end{diagram}$$
equals $\phi:X\to B$. And this holds, because by naturality of $\eta$, $T(\phi)\circ \eta 
_{X}= \eta_{B}\circ \phi$; but $\beta \circ \eta _{B}$ is the 
identity map on $B$, by the unitary law for the algebra structure 
$\beta$.  The second equation is then immediate.

\medskip

 An 
extensive quantity of the form $\delta_{x}$ has total $1\in T({\bf 1})$, 
\begin{equation}\label{totdelx}\tot (\delta_{x})=1.\end{equation}
For, the 
composite
map
$$\begin{diagram}X&\rTo^{\eta_{X}}&T(X)&\rTo ^{T(!)}&T({\bf 1})=R
\end{diagram}$$
 equals the 
composite
\begin{equation}\label{const1x}\begin{diagram}X&\rTo ^{!}&{\bf 1}&\rTo^{\eta _{{\bf 
1}}}&T({\bf 1})=R,\end{diagram}\end{equation}
 by  naturality  of $\eta$ 
w.r.to $!:X\to {\bf 1}$, and this is the map taking value $1$ for 
all $ x\in X$.  

\begin{prop}\label{tottxx}
For the total of $P\in T(X)$, we have
$$\tot (P)= \langle P,1_{X}\rangle ,$$
where $1_{X}$ denotes the function $X\to R$ with constant value 
$1\in R$.
\end{prop}
{\bf Proof.}  The map $1_{X}$ is displayed in (\ref{const1x}) above,
so $\langle P,1_{X}\rangle $ is by definition the result of applying to 
$P\in T(X)$ the  composite
$$\begin{diagram}T(X)&\rTo ^{T(!)}&T({\bf 1})&\rTo ^{T(\eta_{1})}&T^{2}({\bf 1})&\rTo 
^{\mu_{1}}&T({\bf 1})=R
\end{diagram}.$$
 Now by one of the unit laws for a monad, the composite of the two 
last maps here is the identity map of $T({\bf 1})$, so the displayed 
composite is just $T(!)$; this is the map which to $P$ returns the 
total of $P$.

\medskip

\noindent{\bf Notation.} We attempted to make the notation as standard 
as possible. On three occasions, this 
forces us to have double notation, like the $\eta$-$\delta$ doubling 
above, and later $E$ (``expectation'') for $\mu$, and an ``integral'' 
symbol for the pairing of extensive and intensive quantities. An 
exception to standard notation is that the exponential object $B^{X}$ is 
denoted $X\p B$, to keep it online. (Other online notations have also 
been used, like $[X,B]$ or $X\multimap B$.) 
\section{Tensor product of extensive quantities}\label{TPEQX}

Let ${\mathbb T}$ be an algebraic theory, and let $P$ and $Q$ be an 
$X$-ary and a $Y$-ary operation of it. Then one can define, 
semantically, an $X\imes Y$-ary operation $P\otimes Q$:
 given an $X\times Y$-tuple $\theta$ on the ${\mathbb T}$-algebra $B$; 
we think of $\theta$ as an  matrix of elements of $B$ with $X$ rows 
and $Y$ columns. Now evaluate $P$ on each of the  $Y$ columns; this 
gives a $Y$-tuple of elements of $B$; then evaluate $Q$ on this 
$Y$-tuple; this gives an element of $B$. This element is declared to be 
the value of the $X\times Y$-ary operation $P\otimes Q$ on $\theta$.

One might instead first have evaluated $Q$  on each of the rows, and 
then evaluated $P$ on the resulting $X$-tuple; this would in general 
give a different result, denoted $P\tilde{\otimes}Q$; the theory is 
called commutative if $\otimes$ and $\tilde{\otimes}$ agree.

This kind of tensor product was formulated monad theoretically, and without 
reference to the semantics involving $T$-algebras, in the author's 
1970-1972 papers, so as to be applicable for any strong monad on any 
cartesian closed category $\E$; 
it takes the form of two maps\footnote{denoted in \cite{MSMCC} by  
$\tilde{\psi}_{X,Y}$ and 
$\psi _{X,Y}$, respetively} natural in $X$ and $Y \in \E$,
$$\begin{diagram}T(X)\times T(Y) & \pile{\rTo ^{\otimes}\\ 
\rTo_{\tilde{\otimes}}}&T(X\times Y);
\end{diagram}$$
the monad is called {\em commutative} if these two maps agree, for 
all $X$ and $Y$. Both $\otimes$ and $\tilde{\otimes}$ make the 
functor $T$ a {\em monoidal} functor, in particular they satisfy a 
well known associativity constraint: $(P\otimes Q)\otimes S = 
P\otimes (Q\otimes S)$, modulo the isomorphisms induced by $(X\times 
Y)\times Z \cong X\times (Y\times Z)$, and similarly for 
$\tilde{\otimes}$. 
(The ``nullary'' part that goes along with $\otimes$, is a map ${\bf 
1}\to T({\bf 1})$.)

Among the equations satsified is  $\otimes \circ (\eta _{X}\times 
\eta _{Y})=\eta _{X\times Y}$, which in the notation with $\delta$ 
reads: for $x\in X$ and $y\in Y$,
\begin{equation}\label{nattens}\delta _{x}\otimes \delta _{y}= 
\delta_{(x,y)}.\end{equation}

If $A=(A,\alpha )$ and $C=(C,\gamma )$ are $T$-algebras, it makes 
sense to ask whether a map $A\times X \to C$ is a $T$-homomorphism 
in the first variable, cf.\ \cite{BCCM}; and similarly it makes sense to ask whether a 
map $X\times A \to C$ is a $T$-homomorphism in the second variable. 
 We shall use the term ``$T${\em -linear} map'' as 
synonymous with $T$-homomorphism; 
this allows us to use the term 
{\em $T$-bilinear} for  a map $A\times B \to C$ which is a $T$-homomorphism in each of the 
two input variables separately (where $A,B$, and $C$ are 
$T$-algebras), and similarly ``{\em $T$-linear in the first variable}'', 
etc.

If $f: X\times Y \to C$ is any map into a $T$-algebra $C$, it extends 
uniquely over $\eta_{X}\times Y$ to a map $T(X)\times Y \to C$ which 
is $T$-linear in the first variable. Similarly $f$ extends uniquely 
over $X\times \eta_{Y}$ to a map $X\times T(Y) \to C$  which is 
$T$-linear in the second variable. However, a map $X\times Y \to C$ does not necessarily 
extend to a $T$-bilinear $T(X)\times T(Y) \to C$; for this, one needs 
commutativity of $T$: 

\section{Commutative monads}\label{CMX}

We henceforth consider a commutative monad $T=(T,\eta ,\mu )$ on $\E$.
The reader may have for instance the free-abelian-group monad in 
mind.

From \cite{BCCM}, we know that commutativity of $T$ is equivalent to 
the assertion that  $\otimes :T(X)\times T(Y) \to T(X\times Y)$ is 
$T$-bilinear, for all $X$ and $Y$. (In the non-commutative case, 
$\otimes$ 
will only be $T$-linear in the second variable, and $\tilde{\otimes}$ 
will only be $T$-linear in the first variable.)

Then if $C=(C,\gamma )$ is a $T$-algebra, any map $f: X\times Y \to C$ 
extends uniquely over $\eta _{X}\times \eta _{Y}$ to a $T$-bilinear 
map $T(X)\times T(Y)\to C$. Since $f$ also extends uniquely over 
$\eta_{X\times Y}$ to a 
$T$-linear $T(X\times Y) \to C$, one may deduce that $\otimes : T(X)\times 
T(Y) \to T(X\times Y)$ is in fact a {\em universal} $T$-bilinear map out of 
$T(X)\times T(Y)$.

If $B=(B,\beta )$ is a $T$-algebra, and $X$ is an arbitrary object, $X\p B$ carries a 
canonical ``pointwise'' $T$-algebra structure inherited from $\beta$. In the category of 
sets, this is just the ``coordinatwise'' $T$-algebra structure on 
$\Pi _{X}B$. 

Let $(A,\alpha)$ and $(B,\beta )$ be $T$-algebras. 
If $\E$ has sufficiently many equalizers, there is a 
subobject $A\p _{T}B$ of $A\p B$, which in the set case consists of 
those maps $A\to B$ which happen to be $T$-linear. With $T$  
commutative, the subobject $A\p _{T}B \subseteq A\p B$ is in fact a 
sub-$T$-algebra, cf.\ \cite{CCGBCM}.

If $A=(A,\alpha )$, $B=(B,\beta)$, and $C=(C,\gamma)$ are 
$T$-algebras, a map $A\times B \to C$ is $T$-bilinear iff 
its transpose $A\to B\p C$ is $T$-linear, and factors through the 
subalgebra $B\p _{T}C$.

In \cite{MEQ}, Theorem 1, we prove that for a commutative monad $T$, 
the exponential adjoint of the pairing $T(X)\times (X\p B) \to B$ is 
a $T$-bilinear map $T(X) \to (X\p B)\p _{T}B$; so therefore also, we 
have
\begin{thm}\label{Tbilinx}The pairing $\begin{diagram} T(X)\times (X\p B) & 
\rTo^{\langle -,-\rangle}&B\end{diagram}$ is $T$-bilinear, for any 
$T$-algebra $B$. 
\end{thm}

\section{Convolution}\label{ConvX} If $a:X\times Y \to Z$ is any map, and $P\in T(X)$, $Q\in T(Y)$, we 
may form
$a_{*}(P\otimes Q) \in T(Z)$, called the {\em convolution} of $P$ and 
$Q$ {\em along} $a$, and denoted $P*_{a}Q$. Since $\otimes$ is $T$-bilinear and $a_{*}=T(a)$ is 
$T$-linear, it follows that $P*_{a}Q$ depends in a $T$-bilinear way on 
$P, Q$. In diagram, convolution along $a$ is the composite
$$\begin{diagram}T(X)\times T(Y)&\rTo ^{\otimes}&T(X\times Y)&\rTo 
^{T(a)}&T(Z)\end{diagram}.$$
The convolution along the unique map ${\bf 1}\times {\bf 1}\to {\bf 1}$ gives 
a multiplication on $R=T({\bf 1})$, which is commutative.

A consequence of the naturality of $\otimes$ w.r.to the maps $!:X\to 
{\bf 1}$ and $!:Y\to {\bf 1}$  is that \begin{equation}\label{tottensx}
\tot (P\otimes Q)=\tot (P)\cdot \tot (Q),
\end{equation}
where the dot denotes the product in $R=T({\bf 1})$. 
Note that this product is itself (modulo the 
identification $T({\bf 1}\times {\bf 1}) \cong T({\bf {\bf 1}})$) a tensor 
product, $T({\bf 1})\times T({\bf 1})\to T({\bf 1}\times {\bf 1})\cong T({\bf 1})$.

We have, for $x\in X$ and $y\in Y$ 
\begin{equation}\label{natconvx}\delta _{x}*_{a}\delta_{y}= \delta 
_{a(x,y)}.\end{equation}
This follows from (\ref{nattens}), together with naturality of $\eta$ 
w.r.to $a$, which in $\delta$-terms reads
$$T(a)(\delta _{(x,y)})= \delta _{a(x,y)}.$$
If $a$ is an associative operation $X\imes X \to X$, it follows from 
properties of monoidal functors that the convolution along $a$, $T(X)\times T(X) 
\to T(X)$ is likewise associative. If $a$ is commutative, 
commutativity of convolution along $a$ will be a consequence, 
but here, one uses the assumption that the monad $T$ is  commutative.

\section{The space of scalars $R:= T({\bf 1})$}\label{SOSX}

 The space $T({\bf 1})$ plays the role of the ``ring'' of scalars, 
or number line. 
It has a $T$-linear structure, since it is a $T$-algebra, and it 
carries a $T$-bilinear multiplication $m$, namely
$$\begin{diagram}T({\bf 1})\times T({\bf 1})&\rTo^{\otimes}& T({\bf 1}\times {\bf 1}) \cong 
T({\bf 1})\end{diagram},$$
which is commutative and associative. 
The 
multiplicative unit is picked out by $\eta_{{\bf 1}}:{\bf 1}\to T({\bf 1})=R$. - Also $R$ 
acts on any $T(X)$, by
\begin{equation}\label{actxx}\begin{diagram}T({\bf 1})
\times T(X)&\rTo^{\otimes}& T({\bf 1}\times X) \cong 
T(X)\end{diagram},\end{equation}
generalizing the description of multiplication on $R$; it is denoted 
just by a dot $\cdot$. This action is 
likewise $T$-bilinear, and unitary and associative.
These assertions follow from $T$-bilinearity of $\otimes$, and the 
compatibility of $\otimes$ with $\eta$.

There are {\em properties} of 
$T$ which will imply that $T$-algebras carry abelian group  structure; 
in this case, $R$ is a commutative ring, with the above $m$ as 
multiplication, and any $T(X)$ is an $R$-module, with $T$-linear maps 
being in particular $R$-linear, see Section \ref{ASX}. 

\section{Intensive quantities,  and their action on extensive 
quantities}\label{IQAEQX}The space $R=T({\bf 1})$ is a $T$-algebra, with a commutative 
$T$-bilinear commutative monoid structure (so in the additive case, 
Section \ref{ASX}, 
it is in particular a commutative ring). From general 
principles\footnote{it is an aspect of the strength of the monad $T$ 
that the category $\E^{T}$ of $T$-algebras is $\E$-enriched;  
it is even {\em cotensored} (cf.\ \cite{Borceux}) over $\E$; and then $X\p B$ is the 
cotensor of the space $X$ with the $T$-algebra $B$, for any 
$T$-algebra $B$. This then in particular applies to $B=T({\bf 1})$. 
For details, see e.g.\ \cite{MEQ}.} 
follows that for any $X$, the space $X\p R = X\p T({\bf 1})$ inherits a 
$T$-algebra structure and a
$T$-bilinear monoid structure; in fact $-\p R$ is a contravariant 
functor with values in the category of monoids whose multiplication 
is $T$-bilinear. If $f:X\to Y$, the map $f\p R : Y\p R \to 
X\p R$ preserves this structure. The map $f\p R$ is denoted $f^{*}$, and is a 
kind of companion to the covariant $f_{*}:T(X)\to T(Y)$. In the 
terminology of Lawvere \cite{CSQ}, $X\p R$ is a space of {\em 
intensive} quantities on $X$.
Note that $T$ only enters in the form of $R=T({\bf 1})$.

The monoid $X\p R =X\p T({\bf 1})$ acts on any space of the form $X\p T(Y)$, by a 
simple ``pointwise'' lifting of the action of $T({\bf 1})$ on $T(Y)$, 
described in (\ref{actxx}), (with $Y$ instead of $X$):
$$\begin{diagram}\{ X\pitchfork T({\bf 1})\} \times \{X\pitchfork T(Y)\}\cong 
X\pitchfork \{T({\bf 1})\times 
T(Y)\}&\rTo^{X\p \{\otimes_{{\bf 1},Y}\}}&X\p T({\bf 1}\times Y)&\rTo^{\cong}&X\p T(Y)
\end{diagram}.$$
(The monoid structure on $R=T({\bf 1})$ is a speial case.)

\medskip
We shall describe an action of the monoid $X\p R$ on $T(X)$. It has a 
special case the ``multiplication of a distribution by a function'' 
known from (Schwartz) distribution theory\footnote{also, it is 
analogous to the cap-product action of the cohomology ring (cup product) on 
homology, cf.\ e.g.\  \cite{HW}}. Notationally, we let the action be from 
the right, and denote it $\vdash$,
$$\begin{diagram} T(X) \imes (X\p R) & \rTo ^{\vdash }& 
T(X).\end{diagram}$$
It is $T$-linear in the first variable\footnote{it is actually 
$T$-bilinear.}: it is the 1-$T$-linear 
extension   over $\eta _{X \imes R}$ of a certain map 
$X\times (X\p R) \to T(X)$; in other words, 
we describe first $P\vdash \phi$ for the case where $P=\eta (x) = 
\delta _{x}$ for some $x\in X$, and where $\phi :X\to R$ is any function. 
Namely, we put
$$\delta _{x}\vdash \phi := \phi (x)\cdot \delta_{x},$$
recalling that $T(X)$ 
carries a (left) action $\cdot$ by $R=T({\bf 1})$ (of course ``left'' and 
``right'' does not make any difference here, since the monoids in question 
are commutative).

Since $T$-linearity implies homogeneity w.r.to multiplication by 
scalars $\lambda$ in $R$, we have in particular that
\begin{equation}
\label{homogx}
\lambda \cdot (\delta_{x} \vdash \phi )= (\lambda \cdot \delta 
_{x})\vdash \phi
\end{equation}
We shall prove that the action  $\vdash$ is unitary and associative.  
 The unit of $X\p R$ is $1_{X}$, i.e.\  the function with constant 
value $1\in R$. So we should prove $P\vdash 1_{X} =P$. By $T$-linearity in 
the first variable, it  is enough 
to see it when the input $P$ from $T(X)$ is of the form $\delta_{x}$, so 
we should prove, for the unitary property,
$$\delta _{x}\vdash 1_{X} = \delta _{x}$$
and similarly, for the associative property, it suffices to prove
$$(\delta _{x}\vdash \phi )\vdash \psi = \delta_{x}\vdash (\phi \cdot 
\psi )$$
for $\phi$ and $\psi $ arbitrary functions $X\to R$.  The first equation then is a consequence of 
$1_{X}(x)=1$; unravelling similarly the second equation, one sees 
 that the two sides are, respectively
$(\phi (x)\cdot \psi (x))\cdot \delta_{x}$ (using (\ref{homogx})), 
and $(\phi \cdot \psi )(x)$, and their equality is a 
consequence of the pointwise nature (and commutativity) of the multiplication in $R$.

\medskip

The following result serves in the Schwartz theory in essence as the 
{\em definition} of the action $\vdash$ of intensive quantities on 
extensive ones. 
Recall that the multiplicative monoid of $R$ acts on any $T$-algebra 
of the form $T(Y)$, via $\otimes _{1,Y}:T({\bf 1})\times T(Y)\to T({\bf 1}\times 
Y) \cong T(Y)$. The action is temporarily denoted $\dashv$; for 
$Y={\bf 1}$, it is just the multiplication on $R$. This action extends 
pointwise to an action $\dashv$ of $X\p R$ on $X\p T(Y)$.
\begin{prop}\label{assxx} Let $\phi$ be a 
function $X\to R$, and let 
$\psi$ be a function $X\to T(Y)$. Then for any $P\in T(X)$,
$$\langle P\vdash \phi , \psi \rangle = \langle P, \phi \dashv \psi 
\rangle \in T(Y).$$ 
\end{prop}
{\bf Proof.} Since both sides of the claimed equation depend in a 
$T$-linear way on $P$, it suffices to prove the equation for the case 
where $P$ is $\delta_{x}$ for some $x\in X$.
We calculate the left hand side:
$$\langle \delta_{x}\vdash \phi ,\psi \rangle = \langle \phi (x)\cdot 
\delta_{x}, \psi \rangle = \phi (x) \cdot \langle \delta _{x},\psi 
\rangle = \phi (x) \cdot \psi (x),$$
using Proposition \ref{pdelx}, and the right hand side similarly calculates
$$\langle \delta_{x},\phi \dashv \psi \rangle = (\phi \dashv \psi 
)(x),$$
which is likewise $\phi (x)\dashv \psi (x)$, because of the pointwise character of the 
action $\dashv$ on $X\p T(Y)$.
\begin{corol}The pairing $\langle P, \phi \rangle$ (for $P\in T(X)$ 
and $\phi \in X\p R$) can be described in terms of $\vdash$ as 
follows:
$$\langle P, \phi \rangle = \tot (P\vdash \phi ).$$
\end{corol}
{\bf Proof.} Take $Y=1$ (so $T(Y)=R$), and take $\psi = 1_{X}$. Then
$$\tot (P \vdash \phi)= \langle P \vdash \phi, 1_{X}\rangle = \langle 
P, \phi\cdot 1_{X}\rangle $$
using (\ref{tottxx}), and then the Proposition. But $\phi\cdot 1_{X} 
= \phi$.

\section{Additive structure}\label{ASX}

We shall in the present Section describe a simple categorical 
property of the monad $T$, which will guarantee that ``$T$-linearity 
implies additivity'', even ``$R$-linearity'' in the sense of a rig $R 
\in \E$ (``rig"= commutative semiring), namely $R=T(1)$. This 
condition will in fact imply that $\E^{T}$ is an additive (or linear) 
category.

We begin with some standard general category theory,
namely a monad $T=(T,\eta , \mu)$ on a category which has finite 
products and finite coproducts. (No distributivity is assumed.)
So $\E$  has an initial 
object $\emptyset$. If $T(\emptyset )\in \E$ is a terminal object,  
then the object $(T(\emptyset ), \mu 
_{\emptyset})$ is a zero object in $\E ^{T}$, i.e.\ it is both 
initial and terminal. It is initial because $T$, as a functor $\E \to 
\E^{T}$, is a left adjoint, hence preserves initials; and since 
$T(\emptyset )={\bf 1}$, it is also terminal (the terminal object in 
$\E^{T}$ being ${\bf 1}\in \E$, equipped with the unique map $T({\bf 
1})\to {\bf 1}$ as 
structure). This zero object in $\E ^{T}$ we denote $0$. Existence 
of a zero object in a category implies that the category has 
distinguished zero maps $0_{A,B}:A\to B$ between any two objects $A$ 
and $B$, namely the unique map $A\to B$ which factors through $0$. 
For $\E^{T}$, we can even talk about the zero map $0_{X,B}:X\to B$, 
where $X\in \E$ and $B=(B, \beta )\in \E^{T}$, namely $0_{X,B}$ is $\eta _{X}$ 
followed by the zero map $0_{T(X),B}:T(X)\to B$. We have a canonical 
map
$X+Y\to T(X)\times T(Y)$: the composite $X \to  X+Y \to T(X)\times 
T(Y)$ is $(\eta _{X},0_{X,T(Y)})$ (here, the first map is the coproduct inclusion 
map ). Similarly, we have a canonical map $Y \to T(X)\times T(Y)$. 
Using the universal property of coproducts, we thus get a canonical 
map $\phi _{X,Y}: X+Y\to 
T(X)\times T(Y)$. It extends uniquely over $\eta _{X+Y}:X+Y \to 
T(X+Y)$ to a $T$-linear map 
$$\Phi _{X,Y}: T(X+Y)\to T(X)\times T(Y),$$
and $\Phi$ is natural in $X$ and in $Y$.
We say that {\em $T:\E \to \E^{T}$  takes binary coproducts to 
products} if $\Phi _{X,Y}$  is an isomorphism (in $\E$ or equivalently in 
$\E^{T}$) for all $X$, $Y$ in $\E$ . Note that the definition presupposed that $T(\emptyset 
)={\bf 1}$; it is the zero object in $\E^{T}$, so that if $T$ takes binary coproducts to products, it in fact takes finite coproducts to 
products, in a similar sense. So we can also use the phrase ``{\em $T$ 
takes finite coproducts to products}'' for this property of $T$. 

We define an ``addition'' map in $\E^{T}$
; it is a map $+: T(X)\times T(X) $ to $T(X)$, namely the composite
$$\begin{diagram}T(X)\times T(X)&\rTo ^{\Phi _{X,Y}^{-1}}& 
T(X+X)&\rTo ^{T(\nabla )}&T(X)\end{diagram}$$
where $\nabla : X+X\to X$ is the codiagonal. 
So in particular, if $in_{i}$ denotes the $i$th inclusion ($i=1,2$) of 
$X$ into $X+X$, we have
\begin{equation}\label{obv}\begin{diagram}
id_{TX} =&\bigl[TX&\rTo ^{T(in_{i})}&T(X+X)&\rTo^{\Phi _{X,X}}&TX\times TX&\rTo 
^{+}&TX\bigr].
\end{diagram}\end{equation}
Note that this addition map is $T$-linear.
Under the identification $T(X)\cong T(X+\emptyset ) \cong T(X)\times 
1$, the equation (\ref{obv}) can also be read: $ T(!):T(\emptyset)\to 
T(X)$ is right unit for $+$, and similarly one gets that it is a left 
unit. 

We leave to the reader the easy proof of associativity and 
commutativity of the map 
$+:T(X)\times T(X) \to T(X)$. It follows that $T(X)$ acquires 
structure of an abelian monoid in $\E ^{T}$ (and also in $\E$).

For an abelian monoid $A$ in any category, we may ask  whether 
$A$ is an abelian {\em group} or not (so there is a ``minus'' 
corresponding to the $+$); existence of such ``minus'' is a {\em property} of 
$A$, not 
an added structure. If $T$ is a monad which takes finite coproducts 
to products, it makes sense to ask whether the canonical monoid 
structure which $T$-algebras in this case have, is actually an abelian 
group structure; it is therefore a property on such $T$, not an added 
structure. We shall henceforth assume this property, since we need 
``minus'' (= {\em difference}) 
 for  {\em differential} calculus.

In \cite{MEQ}, we proved that
\begin{prop}\label{linxx} Every $T$-linear map $T(X)\to T(Y)$ is compatible with 
the abelian group structure.
\end{prop}

We again assume that $T$ is a commutative monad. Recall that we then have the 
$T$-bilinear action $ T(X)\times T({\bf 1}) \to 
T(X)$. It follows from the Proposition that it is additive in 
each variable separately. 

We have in particular the $T$-bilinear commutative multiplication 
$m:T({\bf 1})\times T({\bf 1})\to T({\bf 1})$, likewise bi-additive, 
$m(x+y,z)=m(x,z)+m(y,z)$, or in the notation one also wants to use,
$$(x+y)\cdot z=x\cdot z + y\cdot z,$$
 so that $R=T({\bf 1})$ carries structure of a commutative 
ring. 
We may summarize:

\begin{prop}\label{modulex}Each $T(X)$ is a module over the ring $R=T({\bf 1})$; each 
$T$-linear map $T(X) \to T(Y)$ is an $R$-module morphism.
\end{prop}
It is more generally  true that  $T$-linear maps $A\to B$ (for $A$ and 
$B \in \E ^{T}$) are $R$-module maps. We shall not use this fact.

\medskip

The property of  $T$ that it ``takes finite coproducts to products'' 
accounts for a limited aspect of contravariance for extensive quantities: 
if a space $X$ is a coproduct, $X=X_{1}+X_{2}$, the isomorphism
$T(X_{1}+X_{2})\cong T(X_{1})\times T(X_{2})$ implies that an 
extensive quantity $P\in 
T(X)$ gives rise to a pair of extensive quantities $P_{1}\in T(X_{1})$ 
and $P_{2}\in T(X_{2})$, which one may reasonably may call the {\em 
restrictions} of $P$ to $X_{1}$ and $X_{2}$, respectively. Now, 
restriction is a contravariant construction, and applies as such, for 
intensive quantities, along {\em any} map. For extensive quantities, and for monads 
$T$ of the special kind studied here, it applies only to quite special 
maps, namely to the inclusion maps into finite coproducts, like 
$X_{1}\to X_{1}+X_{2}$. 

We leave to the reader to philosophize over the extent to which, 
given a distribution of smoke in a given room, it makes sense to talk 
about ``the distribution of this quantity of smoke, restricted to the 
lower half of the room''. 

(For distributions in the Schwartz sense, one may construct some 
further ``restriction'' constructions (restriction to open subsets); this is an aspect of the fact 
that the corresponding intensive quantities (= smooth functions) 
admit an ``extension'' construction from closed subsets.)

\medskip

For $u\in R$, we have the translation map $\alpha ^{u}:R\to R$ given 
by $x\mapsto x+u$. If $P\in T(R)$, we have thus also 
$\alpha^{u}_{*}(P)\in T(R)$. 

We have the following reformulation of the translation maps in terms 
of convolution along the addition map $+ :R\times R \to R$:

\begin{prop}For any $P\in T(R)$ and $a\in R$,
$$\alpha^{a}_{*}(P) = \delta_{a}*P = P * \delta_{a}.$$
\end{prop}
{\bf Proof.} The second equation follows from commutativity of $+$. 
To see the first equation, we observe that both sides of $\alpha^{a}_{*}(P) = 
\delta_{a}*P$ depend $T$-linearly on $P$; so it suffices to prove this 
equation for the case where $P$ is of the form $\delta_{b}$ for  
$b\in R$. But $\alpha ^{a}_{*}(\delta _{b})= \delta _{a+b} = \delta 
_{b+a}$.

\medskip

In particular, we see that $\delta_{0}$ is a neutral element for 
convolution, $ \delta_{0}*P = P * \delta_{0}=P$.

\section{Differential  calculus of extensive quantities 
on $R$}\label{DCEQX}

We attempt in this Section to show how some differential 
calculus of extensive quantities $\in T(R)$ may be developed on equal 
footing with the 
standard differential calculus of intensive quantities (meaning here: 
functions defined on $R$). For this, we assume that the monad $T$ on 
$\E$ has the properties described in Section \ref{ASX}, so in 
particular, $R$ is a commutative ring. To have some differential 
calculus going for such $R$, one needs some further assumption.

Consider  a commutative ring $R$. Assume $D\subseteq R$ is a 
subset  satisfying the following ``KL''-axiom: 

\medskip

{\em for any $f :R \to R$, 
there exists a unique $f':R \to R$ such that for all $x\in R$
\begin{equation}\label{klx}f(x+d) = f(x) + d\cdot f'(x)\mbox{\quad 
for all $d\in D$}.\end{equation}}

\medskip 

Example: 1) models of synthetic differential geometry, with $D$ the set of $d\in R$ with 
$d^{2}=0$ (the simple ``Kock-Lawvere'' axiom says (cf.\ e.g.\ 
\cite{SDG}) a little more than 
this, namely it also asks that any function $f:D\to R$ {\em extends} to a 
function $f:R\to R$.) 

2) Any commutative ring, with $D=\{d\}$ for one single 
invertible $d\in R$. In this case, 
for given $f$, the $f'$ asserted by the axiom is the function
\begin{equation*}f'(x)= \frac{1}{d}\cdot (f(x+d)-f(x)),\end{equation*}
the standard difference quotient.

Similarly, if $V$ is an $R$-module, we say that it satisfies KL, if 
for any $f:R\to V$, there exists a unique $f':R\to V$ such that 
(\ref{klx}) holds for all $x\in R$.

In either case, we may call $f'$ the {\em derivative} of $f$.

It is easy to see that {\em any} commutative ring $R$ is a model, 
using $\{d\}$ as $D$, as in Example 2) (and then also, any $R$-module $V$ 
satisfies then the axiom); this leads to some calculus 
of finite differences. Also, it is true that if $\E$ is the category 
of abstract sets, there are {\em no} non-trivial models of the type in Example 1); 
but, on the other hand, there are other cartesian closed categories $\E$ 
(e.g.\ certain toposes containing the category of smooth manifolds, 
cf.\ e.g.\ \cite{SDG}), and where  
a rather full fledged differential calculus for intensive quantities 
emerges from the KL-axiom.

We assume that $R=T(1)$ satisfies the KL-axiom  (for 
some fixed $D\subseteq R$), and also that any $R$-module of the form 
$T(X)$ does so.

\begin{prop}[Cancelling universally quantified $d$s]\label{cuqdx} If $V$ is an 
$R$-module which satisfies KL, and $v\in V$ has the property that 
$d\cdot v =0$ for all $d\in D$, then $v=0$.
\end{prop}
{\bf Proof.} Consider the function $f:R\to V$ given by $t\mapsto 
t\cdot v$.
Then for all $x\in R$ and $d\in D$
$$f(x+d)= (x+d)\cdot v = x\cdot v + d\cdot v,$$ 
so that the constant function with value $v$ will serve as $f'$. On 
the other hand, $d\cdot v=d\cdot 0$ by assumption, so that the equation may be continued,
$$= x\cdot  v +d\cdot 0 $$
so that the constant function with value $0\in V$ will likewise serve 
as $f'$. From the uniqueness of $f'$, as requested by the axiom, then 
follows that $v=0$.

\medskip

 We are 
now going to provide a notion of {\em derivative} $P'$ for any $P\in 
T(R)$. Unlike differentiation of distributions in the sense of Schwartz, 
which is defined in terms of differentiation of test functions $\phi$, 
 our construction does not 
mention test functions, and the Schwartz definition 
$\langle P', \phi \rangle := -\langle P, \phi '\rangle$ comes in 
our treatment out as a {\em result}, see Proposition \ref{switch} below.

For $u=0$,  $P -\alpha^{u}_{*}(P) =0 
\in T(R)$. 
 Assuming that the $R$-module $T(R)$ is KL, we therefore have  for 
any $P\in T(R)$ that there exists a unique $P' \in T(R)$ such that for 
all $d\in D$,
$$d\cdot P' =  P -\alpha^{d}_{*}(P).$$
Since $d\cdot P'$ has total 0 for all $d\in D$, it follows that $P'$ 
has total 0. 

  Differentiation is translation-invariant:
using
$$\alpha ^{t}\circ \alpha ^{s}= \alpha ^{t+s}=\alpha ^{s}\circ \alpha 
^{t},$$
it is easy to deduce that
\begin{equation}\label{transdif}
(\alpha ^{t}_{*}(P))' = (\alpha ^{t})_{*}(P').
\end{equation}
\begin{prop}Differentiation of 
extensive quantities on $R$ is a $T$-linear process.
\end{prop}
{\bf Proof.} Let temporarily $\Delta : T(R)\to T(R)$ denote the differentiation 
process. Consider a fixed $d\in D$. Then for any $P\in T(R)$,
$d\cdot \Delta (P) =d\cdot P'$ is $P-\alpha^{d}_{*}P$; it is a difference of the two 
$T$-linear maps, namely the identity map on $T(R)$ and 
$\alpha^{d}_{*} =T(\alpha ^{d})$, and as such is $T$-linear. Thus for 
each $d\in D$, the map $d\cdot \Delta : T(R)\to T(R)$ is $T$-linear. 
Now to prove $T$-linearity of $\Delta$ means, by monad theory, to prove equality of two 
maps $T^{2}(R)\to T(R)$; and since $d\cdot \Delta$ is $T$-linear,as 
we proved, it follows that the two desired maps $T^{2}(R)\to T(R)$ 
become equal when post-composed with the map ``multiplication by $d$'': 
$T(R)\to T(R)$. Since $d\in D$ was arbitrary, it follows from KL axiom 
for the $R$-module $T(X)$ that the two desired maps are equal, 
proving $T$-linearity.

\medskip
The structure map $T(R)\to R$  of the $T$-algebra $R=T({\bf 1})$ is $\mu 
_{1}:T^{2}({\bf 1})\to T({\bf 1})$. Just as 
$\eta$ plays a special role, with $\eta (x)$ being the Dirac 
distribution $\delta _{x}$, the structure maps for $T$-algebras play 
a role that sometimes deserves an alternative notation and name; thus 
in particular $\mu _{1}:T(R)\to R$ plays in the context of 
probability distributions the role of {\em expectation}, see \cite{MEQ}, and 
we shall here again allow ourselves a doubling of notation and terminology:
$$E(P):= \mu _{{\bf 1}}(P),$$
the {\em expectation} of $P\in T(R)$. It is a scalar $\in R$.

Note that for $a\in R$, 
\begin{equation}\label{expdelx}
E(\delta_{a})=a;
\end{equation}
since $\delta _{a}$ is $\eta _{R}(a)= \eta _{T({\bf 1})}(a)$, and $E = 
\mu_{{\bf 1}}$, this is a consequence of the monad law that
$\mu_{X} \circ \eta_{T(X)}$ is the identity map of $T(X)$ for any 
$X$, in particular for $X={\bf 1}$.

\medskip

\begin{prop}\label{totexp}Let $P\in T(R)$. Then $$  E(P')=-\tot (P).$$
\end{prop}
{\bf Proof.} The Proposition say that  two maps $T(R)\to R$ agree, 
namely  $E\circ \Delta$ and $-\tot $, where $\Delta$, as above, is the 
differentiation process $P\mapsto P'$. Both these maps are 
$T$-linear, so it suffices to prove that the equation holds for the 
case  $P=\delta_{x}$, so we should prove
$$ E(\delta_{x}')=-\tot (\delta _{x}).$$
 By the principle of cancelling 
universally quantified $d$s (Proposition \ref{cuqdx}), it suffices
to prove that for all $d\in D$ that
$$  d\cdot E(\delta_{x}')=- d\cdot \tot (\delta _{x}).$$
The right hand side is $-d$, by 
(\ref{totdelx}).  
The left hand side is 
\begin{align*}E(d\cdot \delta _{x}') &= E(\delta_{x}- \alpha^{d}_{*}\delta_{x})\\
&= 
E(\delta_{x}- \delta_{x+d})\\
&= E(\delta _{x})-E(\delta _{x+d}) = x - (x+d) =-d,\end{align*}
by (\ref{expdelx}). This proves the Proposition.

\medskip

The differentiation process for functions, as a map $R\p V \to R\p 
V$, is likewise $T$-linear, but this important information cannot be 
used in the same way as we used $T$-linearity of the differentiation
$T(R)\to T(R)$, since, unlike $T(R)$, $R\p V$ (not even $R\p R$)  is not known
 to be freely generated by elementary quantities 
like the $\delta _{x}$s. 

Recall that if $F:V\to W$ is an $R$-linear map between KL modules
\begin{equation}
\label{lindiffx} F\circ \phi' =(F\circ \phi)'
\end{equation}
for any $\phi :R\to V$.

\medskip

One can generalize the differentiation of extensive 
quantities on $R$ to a differentiation of extensive quantities on any 
space $X$ equipped with a vector field. The case made explicit is where the 
vector field is $(x,d)\mapsto x+d$  (or $\tfrac{\partial}{\partial 
x}$) 
on $R$.

\begin{prop}Let $P\in T(R)$ and $Q\in T(R)$. Then
$$(P*Q)'= P'*Q= P*Q'.$$
\end{prop}
{\bf Proof.} By commutativity of convolution, it suffices to prove 
that $(P*Q)'= P'*Q$. Both sides depend in a $T$-bilinear way on $P$ 
and $Q$, so it suffices to see the validity for the case where 
$P=\delta_{a}$ and $Q=\delta _{b}$. To prove 
$(\delta_{a}*\delta_{b})'=\delta_{a}'*\delta_{b}$, it suffices to 
prove that for all $d\in D$, 
$$d\cdot (\delta_{a}*\delta_{b})'=d\cdot \delta_{a}'*\delta_{b},$$
and both sides comes out as
$\delta _{a+b}-\delta _{a+b+d}$, using that $*$ is $R$-bilinear.

\subsection*{Primitives of extensive quantities}\label{PEQX}

We noted already in Section \ref{MATAX} that $P$ and $f_{*}(P)$ have 
same total, for any $P\in T(X)$ and $f:X\to Y$. In particular, for 
$P\in T(R)$ and $d\in D$, $d\cdot P' = P-\alpha^{d}_{*}(P)$ has total 
$0$, so cancelling the universally quantified $d$ we get  that $P'$ has 
total $0$.

A {\em primitive} of an extensive quantity $Q\in T(R)$ is a $P\in 
T(R)$ with $P'=Q$. Since any $P'$ has total 0,
a necessary condition that an extensive quantity $Q\in T(R)$ has a primitive is that $\tot 
(Q)=0$. Recall that primitives, in ordinary 1-variable calculus,  are 
also called 
``indefinite {\em integrals}'', whence the following use of the word 
``integration'':

\medskip
\noindent{\bf Integration Axiom.} {\em Every $Q\in T(R)$ with $\tot (Q)=0$ has a unique 
primitive.}

\medskip
(For contrast: for intensive quantities $\phi $ on $R$ (so $\phi 
:R\to R$ is a function), the standard integration axiom is that 
primitives {\em always} exist, but are {\em not} unique, only up to an additive 
constant.) 

By $R$-linearity of the differentiation process $T(R)\to T(R)$, the 
uniqueness assertion in the Axiom is equivalent to the assertion: 
{\em if $P'=0$, then $P=0$}. (Note that $P'=0$ implies that $P$ is 
invariant under translations $\alpha^{d}_{*}(P)=P$ for all $d\in D$.) 
The reasonableness of this latter 
assertion is a two-stage argument: 1) if $P'=0$, $P$ is invariant 
under {\em arbitary translations} $\alpha ^{u}_{*}(P)=P$. 2) if $P$ 
is invariant under all translations, and has compact support, it must 
be 0. (Implicitly here is: $R$ itself is not compact.)

In standard distribution theory, the Dirac distribution $\delta_{a}$ 
(where $a\in R$) has a primitive, namely the Heaviside ``function''; but 
this ``function'' has not 
compact support - its support is a half line $\subseteq R$.

On the other hand, the integration axiom provides a (unique) primitive 
for a distribution of the form $\delta_{a}-\delta_{b}$, with $a$ and 
$b$ in $R$. This primitive is denoted $[a,b]$, the ``interval'' from 
$a$ to $b$; thus, the defining equation for this interval is
$$[a,b]'=\delta_{a}-\delta_{b}.$$
Note that 
the phrase ``interval from \ldots to \ldots " does not imply that we 
are considering an ordering $\leq$ on $R$ (although ultimately, one 
wants to do so).
\begin{prop} The total of $[a,b]$ is $b-a$. 
\end{prop} 
{\bf Proof.} We have
\begin{align*} \tot ([a,b]) &=-E([a,b]')= -E(\delta _{a}-\delta _{b})\\
\intertext{by Proposition \ref{totexp} and the fact that $[a,b]$ is 
a primitive of $\delta _{a}-\delta _{b}$}
&=-E(\delta_{a})+E(\delta _{b}) = b-a,
\end{align*}
by 
(\ref{expdelx}).

\medskip

 It is of some interest to study the sequence of extensive quantities
$$[-a,a], \quad [-a,a]*[-a,a], \quad  [-a,a]*[-a,a]*[-a,a], \quad \ldots ;$$
they have totals $2a, (2a)^{2},(2a)^{3},\ldots$; in particular, if 
$2a=1$, this is a sequence of probability distributions, approaching 
a Gauss normal distribution (the latter, however, has presently no 
place in our context, since it does not have compact support).

\section{Extensive quantities and Schwartz distributions}\label{EQSDX}
Recall from (\ref{semantx}) that $P\in T(X)$ gives rise to an 
$X$-ary operation on any $T$-algebra $B=(B,\beta )$, via
$\langle P, \phi \rangle := \beta 
(\phi_{*}(P))$, where $\phi \in X\p B$. 
The pairing is thus a map
$$\begin{diagram}T(X) \times (X\p B)&\rTo^{\langle -,- 
\rangle}&B.\end{diagram}$$
which  is 
$T$-bilinear (cf.\ Theorem \ref{Tbilinx}, or \cite{MEQ}). We may take the exponential transpose of the pairing; 
this is then a map $\tau_{X}: T(X) \to (X\p B)\p _{T}B$.

The synthetic rendering 
of Schwartz distribution theory is that  $(X\p R)\p_{T}R$ ``is''  the 
space of Schwartz distributions of compact support: $X\p R$ is the 
``space of test functions'' (not necessarily of compact support), and 
$(X\p R)\p_{T}R$ is the space of $T$-linear functionals $X\p R \to R$ on the space 
of such test functions. (In some well adapted models $\E$ of SDG, and for suitable $T$, this 
 can be proved to be an object whose set of global sections is in 
fact the standard Schwartz distributions of compact support on $X$, if $X$ is 
a smooth manifold; cf.\ \cite{MR} Proposition II.3.6 (Theorem of Que 
and Reyes).)

\medskip

\noindent{\bf Remark.} Consider a commutative ring object in a 
sufficiently cocomplete cartesian closed category 
$\E$. Let $T$ be the (strong) monad which to $X$ associates the free 
$R$-module on $X$. Thus in particular $R=T({\bf 1})$. The monad $X\mapsto 
(X\p R)\p _{T}R$ is in general not a commutative monad, so cannot 
agree with $T$, although in some cases, the monad map $\tau 
:T\rightarrow (-\p R)\p _{T}R$ is monic. In \cite{BET}, it is proved 
that $T(X)$ for a special case (convenient vector spaces) can be 
``carved out'' of $(X\p R)\p_{T}R$ by topological means. Other 
investigations, e.g.\ in \cite{PRSFA}, and a Theorem of Waelbroeck, 
describe a class of spaces $X$ which ``perceive'' $\tau_{X}$ to be an 
isomorphism.

\medskip 

To make contact with classical theory and intuition, we introduce, 
for the third time, a doubling of notation (this one is actually 
quite classical); for $P\in T(X)$ and $\phi \in X\p B$ (where $B$ is 
a $T$ algebra), we write
$$\int _{X}\phi (x)\; dP(x):= \langle P, \phi \rangle \in B,$$
with $x$ a dummy variable ranging over $X$.
Thus  Proposition \ref{tottxx} may be rendered
$$\tot (P) = \int _{X} 1\; dP(x).$$
For $B= T({\bf 1})=R$ and $P\in T(R)=T^{2}({\bf 1})$, the $\int$-notation will help 
to motivate the use of the terminology ``expectation 
of $P$'' for $\mu_{{\bf 1}}(P)$; let $\phi:R\to R$ be the identity map.
Then
$$\langle P,\phi \rangle = \mu_{{\bf 1}}(\phi_{*}(P)) = \mu _{{\bf 1}}(P),$$ 
since $\phi_{*}$ is the identity map of $T(R)$.
On the other hand
$$\int _{X}\phi (x)\; dP(x) = \int _{X}x\; dP(x),$$
since $\phi (x) =x$; this is the standard ``integral'' expression for 
expectation $E(P)$ for a probability distribution $P$ on $R$.

\medskip

The map $\tau _{X}: T(X) \to (X\p R)\p_{T}R$ is not necessarily monic; in 
the case of classical Schwartz distributions, it is monic, which 
allows the classical theory to identify extensive quanties in $T(X)$ 
with elements in $(X\p R)\p _{T}R$, and so one avoids having to 
mention $T(X)$ explicitly; the notion  of extensive quantity on $X$ 
is thus made dependent on the notion of intensive quantity (test function) 
on $X$. It is, however, easy to give examples of $T$s where there are 
not sufficiently many ``test functions'' $X\to R$ (with $R=T(1)$) to 
make $\tau _{X}:T(X) \to (X\p R)\p_{T}R$ injective, whence one 
motivation for the study of  $T$, independent of the introduction of $X\p R$.

The injectivity of $\tau _{X}$ may be expressed: 
``To test equality of $P$ and $Q$ in $T(X)$, it suffices to test, for 
arbitrary functions $\phi :X\to R$, whether $\langle P, \phi \rangle = 
\langle Q, \phi \rangle$'', whence the name {\em test} function.

If $\tau_{X}$ is not injective, there are not sufficiently many such 
test functions $X\to R$, but there are enough, if we allow test 
functions with arbitrary $T$-algebras $B=(B,\beta )$ as their 
codomain.
 We have in 
fact
\begin{prop}For any $X\in \E$, there exists a $T$-algebra $B$ so that
$\tau_{X}:X\to (X\p B)\p_{T}B$ is  monic.
\end{prop}
{\bf Proof.} Take $B=T(X)$; then the map $e: (X\p T(X))\p_{T}T(X)\to 
T(X)$ given 
by ``evaluation at $\eta_{X}:X\to T(X)$'' is left inverse for 
$\tau_{X}$. For, if $P\in T(X)$, then to say $e(\tau _{X}(P))=P$ is 
equivalent to saying
\begin{equation}\label{emergx}\langle P, \eta_{X}\rangle 
=P.\end{equation}
Since the structure map of the $T$-algebra $T(X)$ is $\mu _{X}$,
the definition (\ref{semantx}) of the pairing $\langle P, 
\eta_{X}\rangle$ gives  
$\langle P, 
\eta_{X}\rangle=\mu _{X}(\eta _{*}(P))$ (where $\eta$ here is short for $\eta 
_{X}$). However, $\eta _{*}$ is just another notation for $T(\eta )$, 
and $\mu_{X}\circ T(\eta _{X})$ is the identity map on $T(X)$ by one 
of the monad laws. So we get $P$ back when we apply this map to $P$.

\medskip 
Thus, instead of identifying an extensive quantity on $X$ by its action on 
arbitrary test functions $\phi :X\to R$, we identify it by its action 
on {\em one  single} test function, namely the function $\eta _{X}:X\to T(X)$.

\medskip
I conjecture that $T(X)$ is actually the end $\int _{B\in \E^{T}}(X\p 
B)\p _{T}B$.

\medskip

Here is an important relationship between differentiation of 
extensive quantities on $R$, and of functions $\phi : R\to T(X)$; 
such functions can be differentiated, since $T(X)$ is assumed to be KL  as an 
$R$-module. (In the Schwartz theory, this relationship, with $X={\bf 
1}$, serves as {\em 
definition} of derivative of distributions.)
\begin{prop}\label{switch}For $P\in T(R)$ and $\phi \in R\p T(X)$, one has
$$\langle P' ,\phi \rangle = -\langle P, \phi ' \rangle.$$
\end{prop}
{\bf Proof.} We are comparing two maps $T(R)\times (R\p T(X)) \to 
T(X)$, both of which are $T$-linear in the first variable. Therefore, 
it suffices to prove the equality for the case of $P=\delta _{t}$; in fact, by $R$-bilinearity of the pairing, it suffices to prove that for any $t\in R$ and 
$d\in D$, we have
$$\langle d\cdot (\delta_{t})', \phi \rangle =- \langle \delta _{t}, d\cdot \phi'\rangle.$$ 
The left hand side is $\langle \delta_{t}- 
\alpha^{d}_{*}(\delta_{t}), \phi \rangle$, and using bi-additivity of 
the pairing, this gives
$\phi (t)-( (\alpha^{d})^{*})(\phi )(t) = \phi (t) - \phi (t+d)$, which 
is $-d\cdot \phi '(t)$.

\medskip

Proposition \ref{totexp} can be seen as a special case, with $X={\bf 1}$ 
(thus $T(X)=R$), and with $\phi$ the identity function $R\to R$. 
We use the ``integral'' notation. 
Thus $x$ denotes the identity function on $R$. So
\begin{equation*}E(P')= \int _{R}x\; dP'(x)
=-\int_{R}(x)', dP(x)
= -\int _{R}1\; dP(x),\end{equation*}
the middle equality by  Proposition \ref{switch}. This, however, is 
$-\tot (P)$, by Proposition \ref{tottxx}.

\medskip

The relationship between differentiation of extensive and intensive 
quantities on $R$ expressed in  Proposition \ref{switch} may be given a more 
``compact'' formulation, using (\ref{emergx}).
For then we have, for $P\in T(R)$, that
$$P' =  \langle P', \eta_{R} \rangle = -\langle P, \eta_{R}' \rangle.$$
Thus in particular, knowledge of $\eta_{R}' $ gives knowledge of $P'$ 
for any $P\in T(R)$. It also gives knowledge of $\phi'$ for any 
$\phi:R \to V$, with $V$ a $T$-algebra which is KL module. For, any such $\phi$ 
extends over $\eta _{R}$ to a (unique) $T$-linear  $F:T(R)\to V$, so 
$\phi = F \circ \eta _{R}$; therefore
$$\phi ' = (F\circ \eta_{R})' = F \circ \eta_{R}',$$
using (\ref{lindiffx}).

\medskip

Let us calculate $\eta_{R}':R\to T(R)$ explicitly; we have for any 
$d\in D$ and $x\in R$ (writing $\eta$ for $\eta _{R}$)
$$d\cdot \eta'(x) = \eta (x+d)-\eta (x) =\delta _{x+d}-\delta_{x}= 
-d\cdot (\delta_{x})',$$
so cancelling the universally quantified $d$, we get for any $x\in R$ 
that
$$\eta'(x) =-(\delta _{x})'.$$
The first differentiation refers to differentiation of functions, 
the second to differentiation of distributions; it is tempting to 
write the former with a Newton dot; then we get $\stackrel{\bullet}{\eta}(x) 
=-(\delta _{x})'$.

\medskip

The following depends on the Leibniz rule for differentaiating a 
product of two functions; so this is {\em not} valiud under he 
general assumptions of this Section, but needs the further assumption 
of Example 2, namely thet $D$ consists of $d\in R$ with $d^{2}=0$, 
as in synthetic differential geometry. 
We shall then use ``test function'' technique to prove 
\begin{prop}\label{leibnizx}For any $P\in 
T(R)$ and $\phi \in R\p R$,
$$(P\vdash \phi )' = P'\vdash \phi + P\vdash \phi '.$$
\end{prop}
{\bf Proof.} It suffices to prove that for the ``universal'' test function $\eta 
=\eta _{X} : X\to T(X)$, we have
$$(\langle P\vdash \phi )' , \eta \rangle =\langle P'\vdash \phi, \eta 
\rangle  + \langle P\vdash \phi ', \eta \rangle.$$
We calculate:
\begin{align*}\langle (P\vdash \phi )' , \eta \rangle &= -\langle 
P\vdash \phi , \eta '\rangle \mbox{\quad (by Proposition \ref{switch})}\\
&= -\langle P,\phi \dashv \eta' \rangle \mbox{\quad (by Proposition 
\ref{assxx})}\\
&= -\langle P, (\phi\dashv \eta)'-\phi'\dashv \eta \rangle\\
\intertext{using that Leibniz rule applies to any bilinear pairing, 
like $\dashv$,}
&=-\langle P, (\phi\dashv \eta)'\rangle +\langle P,\phi'\dashv \eta 
\rangle \\
&=\langle P', \phi \dashv \eta \rangle + \langle P, \phi' \dashv \eta 
\rangle\\
\intertext{using Proposition \ref{switch} on the first summand}
&=\langle P'\vdash \phi ,\eta \rangle + \langle P\vdash \phi' ,\eta 
\rangle \\
\intertext{using Proposition \ref{assxx} on each summand}
&=\langle P'\vdash \phi + P\vdash \phi' ,\eta \rangle
\end{align*}
 In other words, (replacing $\dashv$ by $\cdot$), the proof looks 
formally like the one from books on distribution theory, but does not 
depend on ``sufficiently many test functions with values in $R$''.

\begin{verbatim}kock@imf.au.dk
\end{verbatim}

\end{document}